\documentclass[12pt]{amsart}

\usepackage{amsmath,amsthm} 
\usepackage{amssymb} 
\usepackage{graphicx}
\usepackage{hyperref}
\newtheorem{theorem}{Theorem}[section]

\newtheorem{definition}{Definition}

\newtheorem{lemma}[theorem]{Lemma}

\begin{document}

\title[]{Spines of 3-Manifolds as Polyhedra \\with Identified Faces}
\author[]{Sim\'on Isaza}
%\thanks{Partially supported by Colciencias, Colombia, grant: 1118-05-11411.}

\address{Sim\'on Isaza, Escuela de Matem\'aticas, Universidad Nacional de Colombia\\ A.A. 3840 Medell\'{\i}n, Colombia}
\email{psisaza@unal.edu.co, psisaza@gmail.com}

%\date{}
%\thanks{Submitted September 13, 2002.}
\begin{abstract}
In this article we establish the relation between the spines of 3-manifolds
and the polyhedra with identified faces. We do this by showing that the
spines of the closed, connected, orientable 3-manifolds can be presented
through polyhedra with identified faces in a very natural way. We also prove
the equivalence between the special spines and a certain type of polyhedra,
and other related results.
\end{abstract} 
\maketitle

\section{Introduction}

In this article we will consider the polyhedra with identified faces and the
spines of 3-manifolds. Spines have been studied broadly by Matveev, among
other mathematicians, and their definition, as well as the main theorems
concerning them, are available in \cite{Mat2}. Further results on spine
theory can be found in \cite{Mat1}, \cite{MOS}, \cite{MR} and \cite{MP}. On
the other hand, a polyhedron with identified faces is a solid polyhedron
with an even number of polygonal faces, such as a platonic or archimedean
solid, in which every $n$-agonal face is identified with another $n$-agonal
face in a nice way. Among other things, the identification of two faces has
to match up vertices with vertices and edges with edges. Also, the
correspondence of the faces must be frontwards. This means that if two
identified faces are oriented in an induced way by an arbitrary orientation
of the polyhedron, then their identification must be orientation reversing.

Two types of spaces are yielded naturally by these polyhedra. The first one,
the space produced by the polyhedron, consists of the space obtained by
taking the polyhedron and gluing every face frontwards with the face to
which it is identified. The second one, the scar of the polyhedron, is
obtained by taking only the boundary of the polyhedron ($S^{2}$) and
performing the same gluing. Polyhedra with identified faces are usefull in
the study of 3-manifolds because these manifolds can be presented by
polyhedra. More specifically, every closed, connected, orientable 3-manifold
can be obtained as the space produced by some polyhedron, as it is asserted
in Theorem \ref{T1}.

Polyhedra with identified faces are fairly common objects in
three-dimentional topology. The lens spaces, for example, are commonly
presented by this kind of polyhedra (see \cite{Rol}). Other examples are the
dodecahedron which is the base for constructing the Poincar\'{e}'s
dodecahedral space, the fundamental domains of isometry groups of $E^{3}$
and $H^{3}$, and the butterflies, developed and studied by Hilden,
Montesinos, Tejada, and Toro (see \cite{HMTT}, \cite{Tor}). The methods
developed by Montesinos in the study of 3-manifolds as coverings of $S^{3}$
branched over knots have also a close relation with polyhedra with
identified faces (see \cite{Mon1}, \cite{Mon3}). Some other results
concerning this kind of polyhedra can be found in \cite{HLM}\ and \cite{Mon4}%
. However, the study of polyhedra on these works has been mostly lateral,
and has been restricted to particular polyhedra or families of these. In
spite of its broad applications, polyhedra with identified faces have been
rarely defined in a general way or studied independently. This initiative
has been undertaken by Cannon, Floyd, and Parry (see \cite{CFP1}, \cite{CFP2}%
, \cite{CFP3}) and, from a rather different approach, by the author (see 
\cite{Isa}).

The purpose of this article is to establish the relation between the
polyhedra with identified faces and the spines of 3-manifolds. This relation
is very close since, as we will see, the spines of 3-manifolds can be
presented in terms of polyhedra in very convenient ways. This allows us to
use polyhedra to define spines and special spines in ways different from the
classical one, providing tool in the study of spines.

This presentation of spines through polyhedra is useful for several reasons.
In the first place, since any sphere or polyhedron can be seen as a
compactified bidimensional Euclidean space, polyhedra with identified faces
provide bidimensional descriptions (or diagrams) of spines, being spines
intrincated objects of higher dimensions. Besides, complicated mofidications
of spines inside the ambient manifolds can be translated into modifications
of polyhedra which are very easy to perform (see \cite{Isa}). Also, since
polyhedra are closely related to the presentation of manifolds as branched
coverings of $S^{3}$ (see \cite{Mon1}, \cite{Mon3}), this work sets a
connection between spines and such presentations. Finally, since polyhedra
can be frequently embedded on $E^{3}$ and $H^{3}$ nicely (for example when
they are fundamental domains, or when they are realized through Andreev's
Theorem), the results proved here suggest the possibility to endow spines
with some geometrical meaning.

We can mention also that it has been pointed out by Cannon, Floyd, and Parry
that it is very difficult to obtain examples of polyhedra with identified
faces that produce 3-manifolds instead of pseudomanifolds (see \cite{CFP1}, 
\cite{CFP2}). In \cite{Isa} we have already described a method for
constructing virtually an unlimited amount of polyhedra that produce
3-manifolds. Here we also approach the problem by identifying a broad family
of polyhedra that always produce 3-manifolds (Theorem \ref{T2}).

The structure of the article is as follows. In a first section of
preliminary concepts we will introduce briefly the polyhedra with identified
faces. A more complete treatment can be found in \cite{Isa}. In the second
section we will show that the set of homogeneous bidimensional spines of
closed, connected, orientable 3-manifolds, is exactly the set of scars of
the polyhedra that produce such manifolds (Theorem \ref{T2}). Moreover, we will
see that a spine of a 3-manifold $M$ can be seen as the scar of a polyhedron
that produces that same manifold $M$. In the third section we will define a
type of polyhedra that we will call distinguished, and we will prove that
the special spines are exactly the scars of the distinguished polyhedra
(Theorem \ref{T3}). In the fourth section we will see that distinguished
polyhedra are in fact a presentation of the closed, connected, orientable
3-manifolds (Theorem \ref{T5}). In the fifth section we will group the
polyhedra into classes, that we will call alikeness classes, and we will
prove the equivalence between the alikeness classes of polyhedra and the
pairs of the form $(M,S)$, where $M$ is a closed, connected, orientable
3-manifold, and $S$ is a homogeneous bidimensional spine of $M$ (Theorem %
\ref{C1}). Finally, in the sixth section, we will prove the equivalence
between the distinguished polyhedra and the special spines (Theorem \ref{C1}%
), showing that every special spine can be thought as a distinguished
polyhedron and vice versa, and that the presentation of the 3-manifolds by
special spines (or special thickenable PL polyhedra) is equivalent to the
presentation by distinguished polyhedra.

\section{Preliminaries}

We begin by fixing some definitions and notation.

\begin{definition}
\label{D1} Let $K$ be a CW complex of dimension $n$. We say that $K$\ is
homogeneous if every cell of dimension less than $n$ is contained in the
closure of a cell of higher dimension, and if for every cell $\sigma$ of
dimension less than $n-1$, $star(\sigma)-\sigma$ is connected.
\end{definition}

Let us notice that a finite homogeneous CW complex of dimension $n$\ can be
seen in a natural way as a quotient space of $k$\ $n$-dimensional closed
balls, where $k$\ is the number of $n$-dimensional cells in the complex.

\begin{definition}
\label{D2}We say that a topological space $G$ is a graph if it can be
splitted into open cells in such a way that, if $T$ is the set of these
cells, then $\left( G,T\right) $ is a homogeneous one-dimensional CW
complex. In this case we also say that $T$ is a triangulation of $G$. We
call the zero-dimensional and one-dimensional cells of $T$, respectiveley,
the vertices and edges of $G$.
\end{definition}

Let us notice that this definition of graph allows what in graph theory is
known as loops and multiple edges.

\bigskip

\noindent \textbf{Notation.} \label{N1}Given a topological space $X$, we
denote the interior and closure of $X$ by $\overset{\circ }{X}$ and $\bar{X}$
respectiveley. Besides, given an arbitrary natural number $m$, we denote by $%
P_{m}$ the polygon in the complex plane whose vertices are the $m$-th roots
of the unit. We understand $P_{2}$ as the unit disk with two vertices at 1 and -1,
and we understand $P_{1}$ similarily.

\bigskip

We give now a formal definition of a polyhedron with identified faces. Let $%
B^{3}$ be the closed three dimensional ball, $G$ a connected graph embedded
on $\partial B^{3}$, and $T$ a triangulation of $G$. Then we call the terna $%
\left\langle B^{3},G,T\right\rangle $ a \textit{cell-divided ball}. The
following lemma is intuitively clear. A proof can be found in \cite{Isa}.

\begin{lemma}
\label{L1}Let $G$ be as above. Then $\partial B^{3}-G$ is the union of a
finite number of disjoint bidimensional open disks.
\end{lemma}

Let then $\left\langle B^{3},G,T\right\rangle $ be a cell-divided ball such
that $\partial B^{3}-G$ consist of $k$ disjoint open disks. We call the
vertices and edges of $G$ respectiveley the \textit{vertices} and \textit{%
edges} of $\left\langle B^{3},G,T\right\rangle $, and the clousures of the $%
k $ mentioned open disks the \textit{faces} of $\left\langle
B^{3},G,T\right\rangle $.

Let us consider now a cell-divided ball $\left\langle B^{3},G,T\right\rangle 
$ with an even number of faces $k=2n$ oriented in an induced way by an
arbitrary orientation of $B^{3}$. And additionally, for each $m\in 
%TCIMACRO{\U{2115} }%
%BeginExpansion
\mathbb{N}
%EndExpansion
$, let us consider $P_{m}$ with an arbitrary orientation. Now, let us
suppose that the $2n$ faces of $\left\langle B^{3},G,T\right\rangle $ can be
matched by pairs in such a way that, for each pair $\left\{
F_{i},F_{i}^{-1}\right\} $, there exist functions $f_{i}^{+}:P_{m_{i}}%
\longrightarrow F_{i}$ and $f_{i}^{-}:P_{m_{i}}\longrightarrow F_{i}^{-1}$,
where $m_{i}$ is some integer, that satisfy the following conditions:\newline

\begin{itemize}
\item[I.] Both $f_{i}^{+}$ and $f_{i}^{-}$ restricted $\overset{\circ }{P}%
_{m_{i}}$ are homeomorphisms$.$

\item[II.] Both $f_{i}^{+}$ and $f_{i}^{-}$ send vertices to vertices and
edges to edges.

\item[III.] Of $f_{i}^{+}$ and $f_{i}^{-}$, one preserves the orientation
and the other reverses it.

\item[IV.] Both $f_{i}^{+}$ and $f_{i}^{-}$, restricted to any single edge of $%
P_{m_{i}}$, are PL.
\end{itemize}

\bigskip

The fourth condition is licit because the edges of $G$ are, naturally,
homeomorphic to line segments. We can set then fixed homeomorphisms from $(0,1)$ to each edge of $G$ (And also to the two edges of $P_{2}$ and the single edge of $P_{1}$), and think about these edges as PL polyhedra according to these parametrizations. This way, and since $f_{i}^{+}$ and $f_{i}^{-}$ send edges on
edges, we can ask these functions to be PL on the edges of $P_{m_{i}}$.

Now, under these conditions, given a pair of faces $\left\{
F_{i},F_{i}^{-1}\right\} $ of $\left\langle B^{3},G,T\right\rangle $, the
functions $f_{i}^{+}$ y $f_{i}^{-}$ allow us to define a relation $\epsilon
_{i}$ between $F_{i}$ and $F_{i}^{-1}$ under the following rule: $x\in F_{i}$
is related to $y\in F_{i}^{-1}$ if and only if there exists $z\in P_{m_{i}}$
such that $f_{i}^{+}(z)=x$ and $f_{i}^{-}(z)=y$. Taking for each pair of
faces the respective relation constructed in this way we get a set of
relations $\epsilon :=\{\epsilon _{1},...,\epsilon _{n}\},$ that we will
call an \textit{identification scheme,} for the cell-divided ball $%
\left\langle B^{3},G,T\right\rangle $. The notation $\epsilon =\{\epsilon
_{1},...,\epsilon _{n}\}$ will be recurrent. Let us notice that every $%
\epsilon _{i}$ is one to one (1-1) on $F_{i}-G$, and at most 2-1 or 1-2 on
the edges of $G$. In fact the relation $\epsilon _{i}$ identifies
homeomorphically the interior of $F_{i}$ with that of $F_{i}^{-1}$. It also
identifies homeomorphically the edges of $F_{i}$ with edges of $F_{i}^{-1}$.
However, two edges of $F_{i}$ can eventually be identified with a single
edge of $F_{i}^{-1}$, and vice versa.

If $\left\langle B^{3},G,T\right\rangle $ is a cell divided ball and $%
\epsilon $ is an identification scheme for $\left\langle
B^{3},G,T\right\rangle $, we call the quadruple $\left\langle
B^{3},G,T,\epsilon \right\rangle $ a \textit{polyhedron with identified faces%
}. Throughout this article we will use the word polyhedron to refer to a
polyhedron with identified faces, unless otherwise specified. It is worth to
notice that there are cell-divided balls for which there does not exist any
identification scheme, for example the balls with an odd number of faces.
Now, given a polyhedron $\left\langle B^{3},G,T,\epsilon \right\rangle $,
the orbits of the points of $\partial B^{3}$ under $\epsilon $ induce an
equivalence relation on $\partial B^{3}$, that we denote by $Eq\epsilon $
and call the \textit{relation produced by the polyhedron} ($Eq\epsilon $ is
the smallest equivalence relation that contains the relation $\epsilon
_{1}\cup \cdots \cup \epsilon _{n}$). Also, we call the quotient
space $B^{3}\diagup Eq\epsilon =B^{3}\diagup \epsilon _{1},...,\epsilon _{n}$
the \textit{space produced by the polyhedron}.

On the other hand, we have that the relations $\epsilon _{i}$ identify
homeomorphically edges with edges. In this way it is possible to take, on the
set of edges of $\left\langle B^{3},G,T,\epsilon \right\rangle $, the orbit
of a determined edge $a$ under the set of relations $\epsilon $. These orbits are equivalence classes on the
set of edges of the polyhedron. Also, this orbits are cyclic, in the sense that $a_{1},...,a_{m}$ form an orbit
if and only if there are relations $\epsilon _{i_{1}},...,\epsilon _{im}$ such that $%
\epsilon _{i_{1}}$ identifies $a_{1}$ with $a_{2}$, $\epsilon _{i_{2}}$
identifies $a_{2}$ with $a_{3}$, ... , and $\epsilon _{i_{m}}$ identifies $%
a_{m}$ with $a_{1}$. Furthermore, there are no other relations on $\epsilon$ other than $\epsilon _{i_{1}},...,\epsilon _{im}$ identifying edges on $\left\{ a_{1},...,a_{m}\right\} $.

Similarly, since relations $\epsilon _{i}$ identify vertices with vertices, the orbits of these vertices 
under $\epsilon $ induce equivlence classes on the set of vertices. Then,
given an edge $a$ and a vertex $v$, we call the cardinal of the class of $a$
the \textit{cycle} of $a$, and the cardinal of the class of $v$ the \textit{order} of $v$.

At this point we have already a precise definition of a polyhedron with
identified faces. However, to ease the later work we shall narrow the
definition slightly further, as we will do now.

\begin{definition}
\label{D3}Let $X$ be a topological space, and let $r$ and $q$ be two
equivalence relations on $X$. We say that $r$ and $q$ are alike if there
exists an homeomorphism $f:X\longrightarrow X$ \ such that for every $x$ and 
$y$ in $X$, $(x,y)\in r$ if and only if $(f(x),f(y))\in q$. It is easy to
see that if $r$ and $q$ are alike, then $f^{\prime }:A\diagup
r\longrightarrow A\diagup q$ defined by $f^{\prime }(\left\vert x\right\vert
)=\left\vert f(x)\right\vert $ is an homeomorphism.
\end{definition}

Let $\left\langle B^{3},G,T_{1},\epsilon \right\rangle $ and $\left\langle
B^{3},H,T_{2},\eta \right\rangle $ be two polyhedra with identified faces.
Let us suppose that $Eq\epsilon $ and $Eq\eta $ are alike, and that there
exists an homeomorphism $f:\partial B^{3}\longrightarrow \partial B^{3},$ as
in definition \ref{D3}, such that $f(G)=H$. Then on this case we say that $%
\left\langle B^{3},G,T_{1},\epsilon \right\rangle $ and $\left\langle
B^{3},H,T_{2},\eta \right\rangle $ are \textit{essentially equal} polyhedra.
Clearly, essentially equal polyhedra form equivalence classes on the set of
polyhedra. Besides, since $f$ is an homeomorphism that ``preserves" both the
cell-divided structure and the identification scheme, we can assume simply
that $G=H$ and $\epsilon =\eta $. Hence, two essentially equal polyhedra are
always of the form $\left\langle B^{3},G,T_{1},\epsilon \right\rangle $ and $%
\left\langle B^{3},G,T_{2},\epsilon \right\rangle $.

Let us consider now the following example. Let $\left\langle
B^{3},G,T,\epsilon \right\rangle $ be a polyhedron, and let $\left\{
a_{1},...,a_{n}\right\} $ be a class or orbit of edges on $\left\langle
B^{3},G,T,\epsilon \right\rangle $. Then we can insert at the middle point
of each $a_{i}$ a vertex $v_{i}$, obtaining thus a new triangulation $T_{2}$
of $G$, and a new polyhedron $\left\langle B^{3},G,T_{2},\epsilon
\right\rangle $. Let us notice then that $\left\langle
B^{3},G,T_{1},\epsilon \right\rangle $ and $\left\langle
B^{3},G,T_{2},\epsilon \right\rangle $ are equal in their topological
aspects, and differ only because of some redundant vertices. This example
insinuates, correctly, that even though triangulations of graphs play an
important role in the definition of polyhedra, by determining the edges and
vertices of the same, such triangulations are superfluous from a topological
point of view. The fact is, as we will see, that essentially equal polyhedra
are topologically identical, for which it is desirable to work only with simple class representatives.

We will see now the construction of a standard representative for each class
of essentially equal polyhedra. Let $G$ be a graph and $T$ be a
triangulation of $G$. Then, we say that a vertex of $(G,T)$ is \textit{%
topologically superfluous} if it is adyacent to exactly two edges. Let us
consider now a polyhedron $\left\langle B^{3},G,T,\epsilon \right\rangle $,
and a vertex $v_{1}$ of $\left\langle B^{3},G,T\right\rangle $. Let us
recall that every point in the orbit of $v_{1}$ is also a vertex. We say
that $v_{1}$ is a \textit{needless} vertex in $\left\langle
B^{3},G,T,\epsilon \right\rangle $\ if the following conditions are
satisfied:

\begin{itemize}
\item[(1)] Every point $v_{1},...,v_{n}$ in the orbit of $v_{1}$\
is a topologically superfluous vertex of $G$.

\item[(2)] There does not exist an edge of $T$ whose endings belong both to
the set $\left\{ v_{1},...,v_{n}\right\} $.
\end{itemize}

It is worth to mention, even when we are not interested in proving it, that
(1) implies (2) as long as the space produced by $\left\langle
B^{3},G,T,\epsilon \right\rangle $ is not a lens.

Now, let $T_{0}$ be the triangulation obtained by the removal of all the
needless vertices of $T$. We call $T_{0}$ the \textit{standard triangulation}
of $\left\langle B^{3},G,T,\epsilon \right\rangle $. The good definition of $%
T_{0}$ follows inmediately from (1) and (2). Also, we understand that by the
elimination of each vertex, the two edges adyacent to it are merged into a
single edge.

It can be proved that if $T_{0}$ is the standard triangulation for $%
\left\langle B^{3},G,T,\epsilon \right\rangle $, then $\left\langle
B^{3},G,T_{0},\epsilon \right\rangle $ is a well defined polyhedron. It is
clearly seen thus that $\left\langle B^{3},G,T_{0},\epsilon \right\rangle $
is essentially equal to $\left\langle B^{3},G,T,\epsilon \right\rangle $. We
call $\left\langle B^{3},G,T_{0},\epsilon \right\rangle $ the \textit{%
standard polyhedron} for $\left\langle B^{3},G,T,\epsilon \right\rangle $.
It can be seen also that essentially equal polyhedra have the same standard
polyhedron, for which we will also say that $\left\langle
B^{3},G,T_{0},\epsilon \right\rangle $ is the \textit{standard polyhedron}
for the class of $\left\langle B^{3},G,T,\epsilon \right\rangle $ (More
exactly, what we have is that if $\left\langle B^{3},G,T_{1},\epsilon
\right\rangle $ and $\left\langle B^{3},G,T_{2},\epsilon \right\rangle $ are
essentially equal polyhedra, and if $\left\langle B^{3},G,T_{1,0},\epsilon
\right\rangle $ and $\left\langle B^{3},G,T_{2,0},\epsilon \right\rangle $
are their standard polyhedra, then there exists an homeomorphism $f:\partial
B^{3}\longrightarrow \partial B^{3}$ as in Definition \ref{D3} such that $%
f(G)=G$, and such that $f\mid _{G}$ is a graph isomorphism)

From now on we will understand essentially equal polyhedra as the same
polyhedron. For that reason, instead of working with the set of all polyhedra, we will work with the set of classes of essentially equal polyhedra, or equivalently, the set of standard polyhedra. We will narrow thus our definition of polyhedron with identified faces to include only standard polyhedra, that is, polyhedra without needless vertices.

We will also use the following notation. Let $\left\langle
B^{3},G,T_{0},\epsilon \right\rangle $ be the standard polyhedron for the
class of a given polyhedron $\left\langle B^{3},G,T,\epsilon \right\rangle $%
. Such a class, and therefore its standard polyhedron, is determined
uniquely by the graph $G$ and the identification scheme $\epsilon $, so we
can now denote $\left\langle B^{3},G,T_{0},\epsilon \right\rangle $ simply
by $\left\langle B^{3},G,\epsilon \right\rangle $. We will adopt this
notation from now on.

We will finish this section by examining what kind of spaces are produced by polyhedra with identified faces.
Let $X$ be a compact, connected, second-countable Hausdorff topological
space that satisfies the following conditions:

\begin{itemize}
\item[(a)] Every point of $X$, except perhaps a finite number, has a
neighborhood homeomorphic to a three dimensional open ball.

\item[(b)] If $x\in X$ has not a neighborhood of this type, it has a
neighborhood whose closure is homeomorphic to the cone of a connected sum of
tori.
\end{itemize}

In this case we say that $X$ is a pseudomanifold of type\textit{\ }$\mathcal{%
P}_{1}$, and it can be seen that every pseudomanifold of this type is
triangulable (see \cite{Isa}). We say then that $X$ is of type\textit{\ }$%
\mathcal{P}_{1}o$ if it is orientable, and of type\textit{\ }$\mathcal{P}%
_{1}n$ if it is not. Those points of a pseudomanifold that do not have a
neighborhood homeomorphic to a ball, if they exist, are called \textit{%
singularities}. We state now a theorem that determines exactly
which are the spaces produced by polyhedra with identified faces. We will
not prove the theorem, whose proof can be found in \cite{Isa}.

\begin{theorem}
\label{T1}The set of the spaces produced by polyhedra whith identified faces
is exactely the set of pseudomanifolds of type $\mathcal{P}_{1}o$. Besides,
if a polyhedron produces a pseudomanifold $P$, and $x$ is a singularity of $P$, then $x$ comes from the identification of vertices in the polyhedron.
\end{theorem}

In particular, since every closed, connected, orientable 3-manifold is
trivially a pseudomanifold of type $\mathcal{P}_{1}o$, we have that for
every closed, connected, orientable 3-manifold there exists a polyhedron
that produces it. Also, no manifold of other kind can be produced by a
polyhedron.

\section{Spines and Scars}

Given a polyhedron $\left\langle B^{3},G,\epsilon \right\rangle $, we call
the space $\partial B^{3}\diagup Eq\epsilon $ the \textit{scar} of the
polyhedron. Let us notice that the scar of a polyhedron is a subspace of the
space produced by the same, for $\partial B^{3}\diagup Eq\epsilon \subseteq
B^{3}\diagup Eq\epsilon $. Besides, the scar of a polyhedron has a natural
homogeneous bidimensional CW complex structure, where each cell is obtained
in the following way. Each pair of faces of the polyhedron yields an open $2$%
-cell or \textit{face} in the scar, produced when the interior of the faces
are glued together. Similarily, each class of $n$ edges of cycle $n$ in the
polyhedron yields an open $1$-cell or \textit{edge} in the scar, produced
when the $n$ edges are glued toguether becoming a single one. And finally,
each class of $m$ vertices of order $m$ yields a $0$-cell or \textit{vertex}
on the scar in the same way. Let us notice that, since many polyhedra can
produce the same scar, the CW complex structure just defined depends on the
choice of one of these polyhedra.

Leaving aside scars for a moment, let us recall that every spine of a
3-manifold is a CW complex (For the definition and main properties of spines
see \cite{Mat2}). Then we say that a spine is bidimensional if it is
homeomorphic to a bidimensional CW complex. Similarily we will say that a
spine is homogneous if it is homeomorphic to a homogeneous CW complex.

The goal of this section is to prove the following theorem, which
establishes the equivalence between the scars of polyhedra that produce
manifolds, and the spines of 3-manifolds that are homogeneous bidimensional
CW complexes.

\begin{theorem}
\label{T2}Let $M$ be a closed, connected, orientable 3-manifold, and $%
S\subseteq M$. Then, $S$ is a homogeneous bidimensional spine of $M$ if and
only if it is the scar of a polyhedron that produces $M$.
\end{theorem}

Let us begin by recalling some definitions and notations regarding collapses
on simplicial complexes. Given a simplicial complex $K$ in $%
%TCIMACRO{\U{211d} }%
%BeginExpansion
\mathbb{R}
%EndExpansion
^{k}$, we will denote the union of the elements of $K$ by $\left\vert
K\right\vert $, and the barycentric subdivision of $K$ by $K^{\prime }$. On
the other hand, let us recall that if $K$ is a simplicial complex, we say
that an $n$-simplex $\sigma ^{n}$ in $K$\ is \textit{principal} if it is not
a face of any other simplex in $K$ but itself. Additionally, we say that an $%
(n-1)$-simplex $\rho ^{n-1}$ in $K$ is a \textit{free face} of $\sigma ^{n}$
if \ it is a face of $\sigma ^{n}$, of itself, and of no other simplex. Now, if $\sigma $ is
principal in $K$ and $\rho $ is a free face of $\sigma $, we say that the
complex $K$ \textit{collapses elementally} to the subcomplex $K-\{\sigma
,\rho \}$, and we write $K\searrow K-\{\sigma ,\rho \}$. If $L$ is a
subcomplex of $K$, and it is possible to obtain $L$ form $K$ by means of a
finite secuence of elementary collapses $K\searrow \cdots \searrow L$, we
say that $K$ \textit{collapses} to $L$, and we write $K\searrow L$ (see \cite%
{Mat2}).

We shall now prepare the hypotheses necessary to state a lemma that is
rather technical, but will be necesarry for a step of the proof of Theorem %
\ref{T2}. Let $K$ be a finite three dimensional simplicial complex in $%
%TCIMACRO{\U{211d} }%
%BeginExpansion
\mathbb{R}
%EndExpansion
^{k}$, and $L$ be a bidimensional subcomplex of $K$ with an even number of
triangles $l_{1},...,l_{2n}$. On the other hand, for each $1\leq i\leq n$,
let $f_{i}:%
%TCIMACRO{\U{211d} }%
%BeginExpansion
\mathbb{R}
%EndExpansion
^{k}\longrightarrow 
%TCIMACRO{\U{211d} }%
%BeginExpansion
\mathbb{R}
%EndExpansion
^{k}$ be an affine transformation such that $f_{i}(l_{i})=l_{i+n}$. We will
have under consideration the spaces $\left\vert K\right\vert \diagup \left\{
f_{i}\right\} _{i=1}^{n}$ and $\left\vert L\right\vert \diagup \left\{
f_{i}\right\} _{i=1}^{n}$.

Let us denote by $\pi $ the projection of $\left\vert K\right\vert $ onto $%
\left\vert K\right\vert \diagup \left\{ f_{i}\right\} $. Then, if $\left\{
T_{i}\right\} $ is the set of simplexes of the second barycentric
subdivision of $K$, $K^{\prime \prime }$ (i.e. $\left\{ T_{i}\right\}
=K^{\prime \prime }$), we have that $\left\{ \pi (T_{i})\right\} $ is a
triangulation of $\left\vert K\right\vert \diagup \left\{ f_{i}\right\} $.
We obtain in this way that $\left\vert K\right\vert \diagup \left\{
f_{i}\right\} $ can be embedded in $%
%TCIMACRO{\U{211d} }%
%BeginExpansion
\mathbb{R}
%EndExpansion
^{m}$, for a large enough $m$ (see \cite{ST}), and that $\left\{ \pi
(T_{i})\right\} $ is in fact a simplicial complex such that $\left\vert
\left\{ \pi (T_{i})\right\} \right\vert =\left\vert K\right\vert \diagup
\left\{ f_{i}\right\} $. Similarily, if $\left\{ S_{i}\right\} =L^{\prime \prime }$, then $\left\{ \pi (S_{i})\right\} $ is
a simplicial complex such that $\left\vert \left\{ \pi (S_{i})\right\}
\right\vert =\left\vert L\right\vert \diagup \left\{ f_{i}\right\} $. We are
now in a position to state the following lemma.

\begin{lemma}
\label{L2}Under the previous hypotheses, if $K\searrow L$ then $\left\{ \pi
(T_{i})\right\} \searrow \left\{ \pi (S_{i})\right\} $ (That is, if $%
\left\vert K\right\vert \searrow \left\vert L\right\vert $, then $\left\vert
K\right\vert \diagup \left\{ f_{i}\right\} \searrow \left\vert L\right\vert
\diagup \left\{ f_{i}\right\} $)
\end{lemma}

\textit{Proof}. We know that if $K\searrow L$ then $K^{\prime \prime
}\searrow L^{\prime \prime }$. Let $T$ be a simplex in $K^{\prime \prime }$,
and $\rho $ a free face of $T$ not belonging to $L^{\prime \prime }$. Under
these circumstances it is enough to proove that if $K^{\prime \prime
}\searrow K^{\prime \prime }-\{T,\rho \}$, then $\left\{ \pi (T_{i})\right\}
\searrow \left\{ \pi (T_{i})\right\} -\{\pi (T),\pi (\rho )\}$. The general
result is obtained from here inductively.

Now, it is clear that if $T$ is principal in $K^{\prime \prime }$, then $\pi
(T)$ is principal in $\left\{ \pi (T_{i})\right\} $. Besides, if $\rho $ is
a free face of $T$, the fact that $\rho \notin L^{\prime \prime }$ implies
that $\pi (\rho )$ is a free face of $\pi (T).$ Therefore, if $K^{\prime
\prime }\searrow K^{\prime \prime }-\{T,\rho \}$, then $\left\{ \pi
(T_{i})\right\} \searrow \left\{ \pi (T_{i})\right\} -\{\pi (T),\pi (\rho
)\} $. $\square $

\bigskip

We will now engage in the proof of Theorem \ref{T2}.

\bigskip

\textit{Proof of Theorem }\ref{T2}. For a pseudomanifold $P$\ of type $%
\mathcal{P}_{1}o$ let us define a spine as a generalization of the concept
of a spine of a 3-manifold: If $S\subseteq P$, and $B\subseteq P$ is an open
three dimensional ball with $B\cap S=\phi $, we say that $S$ is a \textit{%
spine} of $P$ if every singularity of $P$ belongs to $S$, and if $%
P-B\searrow S$. It is clear that if $P$ is a 3-manifold, this definition
coincides with the usual definition of a spine of a 3-manifold. Besides,
since every pseudomanifold $P$\ of type $\mathcal{P}_{1}o$ is triangulable,
every spine of $P$ is a CW complex. We will prove now a generalized version of Theorem \ref{T2} that states the
following:

\bigskip 

\textit{If }$P$\textit{\ is a pseudomanifold of type }$P_{1}o$\textit{, and }%
$S\subseteq P$\textit{, then }$S$\textit{\ is a homogeneous bidimensional
spine of }$P$\textit{\ if and only of it is the scar of a polyhedron that
produces }$P$\textit{.}

\bigskip 

Let us see first that if $S$ is the scar of a polyhedron $\left\langle
B^{3},G,\epsilon \right\rangle $ that produces $P$, then $S$ is an
homogeneous bidimensional spine of $P$. We know from Theorem \ref{T1} that
every singularity of $P$ comes from the identification of vertices in $%
\left\langle B^{3},G,\epsilon \right\rangle $. Since these vertices lie in $%
\partial B^{3}$, all singularities of $P$ are contained in $\partial
B^{3}\diagup \epsilon =S$. Then, to see that the scar $S$ is a spine of $P$
we only need to prove that $P$ without a ball collapses to $S$.

Let $U$ be an open ball contained in $B^{3}$ not touching $\partial B^{3}$.
This means that the closure of $U$ is contained in the interior of $B^{3}$.
Now, let $U_{P}$ be an arbitrary open ball in $P$. Without loss of
generality $U_{P}$ comes from $U$, that is, $U\diagup Eq\epsilon =U_{P}$.
Therefore, $(B^{3}-U)\diagup Eq\epsilon =P-U_{P}$. On the other hand, we
have that $\partial B^{3}\diagup Eq\epsilon =S$. Since $(B^{3}-U)\searrow
\partial B^{3}$\ naturally, the previous lemma allows us to conclude that $%
(B^{3}-U)\diagup Eq\epsilon \searrow \partial B^{3}\diagup Eq\epsilon $, or
in other words, that $P-U_{P}\searrow S$. From this we have that $S$ is a
spine of $P$. The fact that $S$ is bidimensional and homogeneous is derived
from the fact, observed at the begining of this section, that every scar is
a homogeneous bidimensional CW complex.

Let us prove now the other implication of the Theorem. Let us see that if $S$
is a homogeneous bidimensional spine of $P$, then $S$\ is the scar of a
polyhedron that produces $P$. The proof will be constructive. The idea is
simply to cut along the spine to obtain a ball and an identification scheme
that re-pastes the cut. In this way we will exhibit the polyhedron required.
Along this proof we will do well by keeping in mind that a binary relation
on $A$ is a subset of $A\times A$. Then, the signs of inclusion and
difference applied to relations will mean nothing but inclusion and
difference of sets.

Let $M$ be a closed, connected,\ orientable 3-manifold $M$, and let $S$ be a
spine of $M$. We know that $M-S$ is a three dimensional open ball (see 
\cite{Isa}, \cite{Mat2}). It is esay to see that the same is true for a
pseudomanifold, with identical proof. Let $P$ be a pseudomanifold of type $%
\mathcal{P}_{1}o$ and let $S$ be a spine of $P$. Since $S$ is a homogeneous
bidimensional CW complex, we can consider a triangulation $T$ of $S$. Using
this we can also consider a triangulation of $P$, whose (closed) tetrahedra $\Delta
_{1},...,\Delta _{k}$ are the triangles of $T^{\prime \prime }$ extended
radially to the center of the ball $B:=P-S$. Let $\Delta _{1}^{\prime
},...,\Delta _{k}^{\prime }$ be a collection of disjoint closed three
dimensional tetrahedra in $%
%TCIMACRO{\U{211d} }%
%BeginExpansion
\mathbb{R}
%EndExpansion
^{3}$. Then $P$ can be viewed as a quotient space $\cup \Delta _{i}^{\prime
}\diagup \sim $, for a certain equivalence relation $(\sim )$ in the union
$\cup \Delta _{i}^{\prime }$.

On the other hand, let us consider the sets $\Lambda _{1},...,\Lambda _{k}$
defined by $\Lambda _{i}=\Delta _{i}\cap B$. Then, there are subsets $%
\Lambda _{1}^{\prime }\subseteq \Delta _{1}^{\prime }$, ... , $\Lambda
_{k}^{\prime }\subseteq \Delta _{k}^{\prime }$, such that $B$ can be seen as
the quotient space $\cup \Lambda _{i}^{\prime }\diagup \simeq $, where $%
(\simeq )$ is the restriction of $(\sim )$ to the union $\cup \Lambda
_{i}^{\prime }$.

Now, if we consider $B$ embedded in $%
%TCIMACRO{\U{211d} }%
%BeginExpansion
\mathbb{R}
%EndExpansion
^{3}$, and if we consider there its closure $\overline{B}$,\ we see that $%
\overline{B}$\ has a natural triangulation induced by $\Lambda
_{1},...,\Lambda _{k}$. Moreover, the elements of this triangulation can be
thought as $\Delta _{1}^{\prime },...,\Delta _{k}^{\prime }$ after suffering
a certain gluing. More specifically, $\overline{B}$ can be seen as a
quotient space $\overline{B}=\cup \Delta _{i}\diagup \approx $, for a
certain extension $(\approx )$ of $(\simeq )$, such that $(\simeq )\subset
(\approx )\subset (\sim )$.

Let $j$ be a number between $1$ and $k$. Let $A_{j}\subset \overline{B}=\cup \Delta _{i}\diagup \approx $
be the very same $\Delta _{j}$ after being glued to others of its kind, according to $(\approx )$, to form $\overline{B}$. That is, let $A_{j}$ the set of points in $\overline{B}=\cup \Delta _{i}\diagup \approx $ corresponding to a class of $%
(\approx )$ that cointains at least an element of $\Delta _{j}$. Let us
consider now the sets $D_{j}$ defined by $D_{j}=A_{j}\cap \partial \overline{%
B}$. Then $D_{j}$ can be understood as the set of points in $\partial 
\overline{B}$ coming from points in $\Delta _{j}$, or just ``the points of $%
\Delta _{j}$ in the boundary of $\overline{B}$". We see that $\cup \partial
D_{i}$, taking these boundaries in $\partial \overline{B}$, is a graph
imbedded in $\partial \overline{B}$, and that for every $i$, $\overset{\circ 
}{D_{i}}$ is a disk. Hence, $\left\langle \overline{B},\cup \partial
D_{i}\right\rangle $ is a cell-divided ball. Furthermore, $(\sim )-(\approx
) $ is a relation in $\partial \overline{B}$ that can be seen as the
relation produced by an identification scheme $\epsilon $ in the faces of $%
\left\langle \overline{B},\cup \partial D_{i}\right\rangle $ (i.e. $[(\sim
)-(\approx )]=Eq\epsilon $). Since%
\begin{equation*}
\begin{array}{c}
\overline{B}\diagup \lbrack (\sim )-(\approx )]%
\end{array}%
=%
\begin{array}{c}
\lbrack \cup \Delta _{i}\diagup \approx ]\diagup \lbrack (\sim )-(\approx )]%
\end{array}%
=%
\begin{array}{c}
\cup _{i}\Delta _{i}\diagup \sim%
\end{array}%
=%
\begin{array}{c}
P,%
\end{array}%
\end{equation*}%
it follows that $\left\langle \overline{B},\cup \partial D_{i}\right\rangle $
under the relation $(\sim )-(\approx )$ is a polyhedron that produces $P$.
Let us notice that $(\sim )-(\approx )$ is a relation on the boundary of $\overline{B}$, implying that $B\diagup \lbrack (\sim )-(\approx )]=B$. Then, since the identification space $B\diagup \lbrack (\sim )-(\approx )]$ is $%
B=P-S $, clearly $\partial B\diagup \lbrack (\sim )-(\approx )]$ is $S$, for
which $S$ is in fact the scar of the polyhedron. $\square $

\bigskip

\section{Special Spines and Scars}

In the previous section we set the equivalence between the homogeneous
bidimensional spines of 3-manifolds and the scars of polyhedra that produce
manifolds. In this section we will establish sufficient and necesary
conditions over a polyhedron for the scar it produces to be a special spine.
It will be necesary to remember that every scar has a natural CW complex,
whose cells are the faces, edges and vertices defined in the first paragraph
of the previous section.

We shall begin with an analysis of the shape of the neighborhoods for the
points of a given scar, produced by a given fixed polyhedron. In the first
place, if a point lies in the interior of a face of the scar, then a closed
regular neighborhood of it in the same scar will be a closed bidimensional
disk (Figure 1 (a)), because the interior of the face of the scar is
produced just by gluing the interiors of two faces of the polyhedron. In the
second place, if a point $x$\ is in the interior of an edge of the scar, and
that edge is produced by the identification of $n$ edges of cycle $n$ of the
polyhedron, then a closed regular neighborhood of $x$\ will be a set of $n$
half closed disks, glued linearily by their diameters (As in Figure 1 (b1)
and (b2)). Particularly, if $x$ is in the interior of an edge produced by
the identification of $2$ edges of cycle $2$, its regular neighborhood will
be a disk (Figure 1 (b1)). Finally, if $x$ is a vertex of the scar produced
by the identification of $m$ vertices of order $m$ of the polyhedron, a
closed regular neighborhood of $x$\ will consist of a series of closed
circular sectors, where one of the two radii in the boundary of each sector,
or both, are glued linearily with other such radii, in a way that the
centers of all the circular sectors end up glued together at a single point,
that in fact is $x$ (as in Figure 1 (c1) and (c2)). In this case the
circumference archs in the boundary of the circular sectors will form a
graph that can be embedded in a compact, connected, orientable 2-manifold
without boundary, dividing the latter in $m$ open disks. Such 2-manifold
will be that whose cone is the regular neighborhood of $x$ in the
pseudomanifold produced by the polyhedron. In this way, if the polyhedron
produces a 3-manifold, it will be possible to embed that graph in $S^{2}$,
that is, it will be a planar graph (as in Figure 1 (c1)).

\begin{figure}[h!]
\centering
    \includegraphics[scale=1]{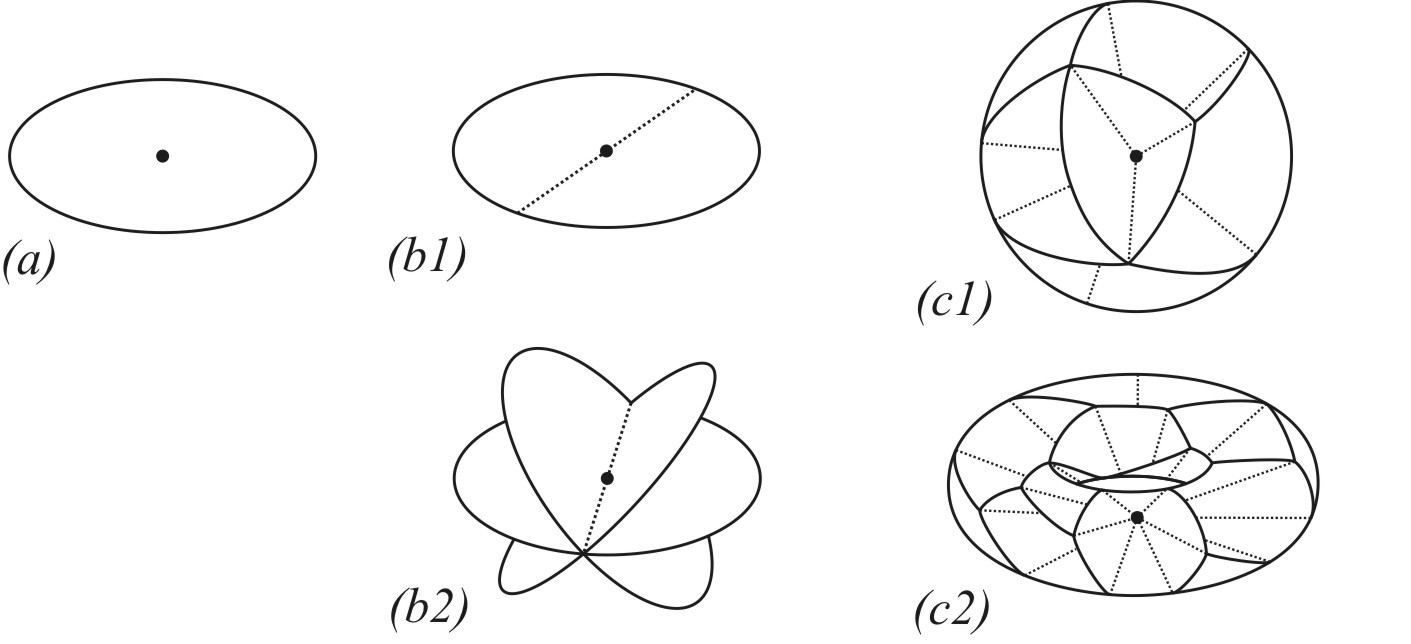}\\
     \caption{Examples of
neighborhoods of points of a scar. \ The neighborhood in (c2) is shown in a
schematic way, for it is\ a certain subset of the conus of a torus that
cannot be embedded in the three dimensional Euclidean space.}
\end{figure}
    
%%%%%%--   

%\FRAME{dtbpFU}{4.7729in}{2.1932in}{0pt}{\Qcb{Figure 1. Examples of
%neighborhoods of points of a scar. \ The neighborhood in (c2) is shown in a
%schematic way, for it is\ a certain subset of the conus of a torus that
%cannot be embedded in the three dimensional Euclidean space.}}{}{%
%espinas_poliedros_0.jpg}{\special{language "Scientific Word";type
%"GRAPHIC";maintain-aspect-ratio TRUE;display "USEDEF";valid_file "F";width
%4.7729in;height 2.1932in;depth 0pt;original-width 4.7184in;original-height
%2.1534in;cropleft "0";croptop "1";cropright "1";cropbottom "0";filename
%'Dibujos/Espinas_poliedros_0.jpg';file-properties "XNPEU";}}

Based on the different shapes of the neighborghoods we define the $2$%
-components, $1$-components and $0$-components of a scar, and in fact of any
bidimensional homogeneous CW complex, in the following way. We define a 
\textit{2-component} of a scar as a connected component of the space of
points whose neighborhoods are disks. Similarly, we define a \textit{%
1-component }as a connected component of the space of points whose
neighborhoods are built from $n$ half disks, as we showed, with $n\neq 2$.
Finally, we define a \textit{0-component} as a set of the form $\left\{
x\right\} $, where $x$ is a point with any other type of neighborhood.

We say then that a scar, or a bidimensional homogeneous CW complex in general, is \textit{cellular} if every $i$-component is an
open cell of dimension $i$. Let us notice that not every scar is cellular.
Let us consider the case of a polyhedron containing a succesion of faces in
which every face limits with the next one along an adge of cycle $2$.
Besides, let us suppose that the last face also limits with the first one
along an edge of cycle $2$, closing a loop. Let us see what happens with the
scar of this polyhedron. Since all the points in the scar coming from the
interior of a face, as well as all the points coming from the interior of an
edge of cycle 2, have neighborhoods with the shape of a disk, it follows
that this ``loop" of faces gives place possibly to a $2$-component with the
shape of an annulus. This example is also useful to show how the $2$%
-components, $1$-components, and $0$-components of a scar do not have to
coincide necessarily with its faces, edges and vertices. 

Now, we say that a scar, or a bidimensional homogeneous CW complex in general $S$, is \textit{simple} if every point in $S$ has a
regular neighborhood (in $S$) shaped like one of the three types of
neighborhoods ilustrated in the following figure. Let us recall now that a spine of a closed connected 3-manifold is called \textit{special} if it is a
homogeneous bidimensional CW complex that is cellular and simple (see \cite{Mat2}).

\begin{figure}[h]
\centering
    \includegraphics[scale=1]{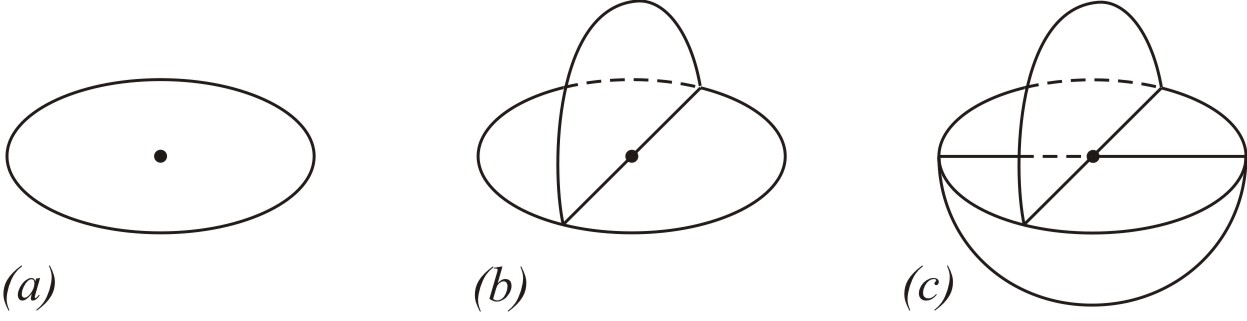}\\
     \caption{ Neighborhoods allowed in special spines.}
\end{figure}

%\FRAME{dtbpFU}{4.2065in}{1.1035in}{0pt}{\Qcb{Figure 2. Neighborhoods allowed
%on special spines.}}{}{espinas_poliedros_1.jpg}{\special{language
%"Scientific Word";type "GRAPHIC";maintain-aspect-ratio TRUE;display
%"USEDEF";valid_file "F";width 4.2065in;height 1.1035in;depth
%0pt;original-width 4.1563in;original-height 1.0706in;cropleft "0";croptop
%"1";cropright "1";cropbottom "0";filename
%'Dibujos/Espinas_poliedros_1.jpg';file-properties "XNPEU";}}

We shall give one more definition with the aim of simplifying the following
proofs. If an edge of a scar $S$ comes from the identification of $n$ edges
of cycle $n$ of the polyhedron, we say that it is an \textit{identified edge
of cycle }$n$. Similarly, if a vertex of $S$ comes from the identidication
of $m$ vertices of order $m$ of the polyhedron, we say that it is an \textit{%
identified vertex of order }$m$. The following four lemmas will be used to
prove subsequent important results, and particulary Theorem \ref{T3} that is
the main result of this section.

\begin{lemma}
\label{L3}If a polyhedron produces a simple scar, then it produces a
3-manifold.
\end{lemma}

\textit{Proof}. Let $P$ be the space produced by the polyhedron, and $S$ its
scar. By Theorem \ref{T1} we know that every singularity of $P$, if there is
any, comes from vertices of the polyhedron and therefore lies in $S$. It
suffices to show then that every point of $S$ has a neighborhood in $P$\
homeomorphic to a ball. For $v\in S$ let $V$ be a closed regular
neighborhood of $v$ in $P$. Since $V$ is the cone of $\partial V$, to show
that $V$ is a ball we only need to prove that $\partial V$ is a sphere.

We know that $V\cap S$ is a closed regular neighborhood of $v$ in $S$. Since 
$S$ is simple, $V\cap S$ has neccesarily one of the shapes (a), (b) o (c)
shown in Figure 2. Besides, $\partial V\cap S$ is a graph embedded in $%
\partial V$ in such a way that $\partial V-S$ consists of open disks. For
this reason, $\partial V\cap S$ induces a CW complex structure in $\partial
V $. If $V\cap S$ has the shape (a), then $\partial V\cap S$ is a
circumference for which $\partial V$ is neccesarily a sphere. If $V\cap S$
has the shape (b), then $\partial V\cap S$ is a graph formed by a
circumference and one diameter, for which $\partial V$ is neccesarily a
sphere. Finally, If $V\cap S$ has the shape (c), then $\partial V\cap S$ is
a complete graph of order 4, that is, a circumference with three radii, for
which, once again, $\partial V$ is a sphere. In every case the conclusion
holds because none of the three graphs can be embedded in a closed,
connected, orientable 2-manifold, other than the sphere, splitting it into
open disks. $\square $

\bigskip

\begin{lemma}
\label{L4}A scar is simple if and only if there exists a polyhedron that
produces it satisfying the following three conditions: All of its edges are
of cycle $2$ or $3$, all of its vertices are of order less than or equal to $%
4$, and all of its vertices are adjacent to exactly three edges of cycle $3$.
\end{lemma}

\textit{Proof}. The implication from left to right follows naturally from
the neighborhood analysis just exposed. Let us see the other implication. We
need to show that given a polyhedron with the three properties stated, then
every point of its scar has a regular neighborhood in that scar with one of
the three shapes allowed. From the previous neighborhood analysis we have
that the points in the faces of the scar have neighborhoods with the shape
of a disk (Figure 2, (a)). We also have that, since the polyhedron has only
edges of cycles 2 and 3, the points in the edges of the scar have
neighborhoods homeomorphic either to a disk or to three half disks glued by
its diameters (Figure 2, (a) y (b)).

It only remains then to examine the case of vertices. Let $v$\ be a vertex
of the scar, then $v$ is by hypotesis an identified vertex of order less
than or equal to 4. However, $v$ does not have order 1 because that would
imply the existence of edges of cycle 1. On the other hand, let us suppose
that $v$ has order 3. Then there exist neccesarily two identified edges of
cycle 3 adjacent to $v$; and $v$ has a neighborhood with the shape of three
half disks, as in  \hbox{Figure 2. (b).} But this implies that the three vertices in
the polyhedron whose identification turns them into $v$ are adjacent to
exactly two edges of cycle 3, which violates the hypotheses. As a
consecuence, $v$ does not have order 3 either.

We conclude that $v$ has neccesarily order 2 or 4. If $v$ has order 2, then
every identified edge in the scar adjacent to $v$ has order 2, and $v$ has a
neighborhood with the shape of a disk.

Let us consider now the case of $v$ having order 4. Let $P$ and $S$ be
respectively the space and scar produced by the polyhedron. Since $P$ is a
pseudomanifold, we know that a closed neighborhood $V$\ of $v$ in $P$ is
homeomorphic to the cone of $\partial V$, which is a closed, connected,
orientable 2-manifold. We know also that $V\cap S$ is a closed regular
neighborhood of $v$ in $S$. Besides, $\partial V\cap S $ is a graph embedded
in $\partial V$ in such a way that $\partial V-S$ consists of open disks,
for which $\partial V\cap S$ induces a CW complex structure in $\partial V$.
Now, since every vertex of the polyhedron is adjacent to exactly three edges
of order 3, we see that $\partial V-S$ consists in fact of open triangles,
and that the CW structure induced in $\partial V$ by $\partial V\cap S$ is a
triangulation. Let us denote by $K$ the triangulation of $\partial V$ thus
obtained.

Since $v$ has order 4, then $K$ has only four triangles. Besides, since
every edge of the polyhedron has cycle 3, every point of $\partial
V=\left\vert K\right\vert $ can lie, at most, in the boundary of three
triangles. The only triangulation of a closed, connected, orientable
2-manifold, made of four triangles, and satisfying this condition, is that
of the sphere triangulated as a tetrahedron. Hence, $V\cap S$ has the shape
shown in Figure 2. (c). $\square $

\bigskip

\begin{lemma}
\label{L5}A scar is cellular if and only of there exists a polyhedron that
produces it not containing edges of cycle $2$.
\end{lemma}

\textit{Proof}. Let $S$ be a cellular scar. Then, since every $i$-component
of $S$ is an open $i$-dimensional cell, the set $C$\ of $i$-components of $S$
endows $S$ with a CW complex structure, or in other words, $(S,C)$ is a CW
complex.

Let us recall that in the proof of Theorem \ref{T2}, to prove that every
spine $\check{S}$\ is the scar of some polyhedron, we started from the fact
that $\check{S}$ had some homogeneous CW complex structure. From there we
proceeded to triangulate $\check{S}$, and then to construct a polyhedron.
However, let us notice that if $S$ is a simple scar, taking $S$ with the CW
complex structure $(S,C)$ that we defined, we can carry out just the same
construction of the proof of Theorem \ref{T2}, but abstaining from
triangulating $S$ (The $\Delta $'s appearing in the proof will not be then
tetrahedra but pyramids with polygonal bases). In this way we obtain a
polyhedron $\left\langle B^{3},G,\epsilon \right\rangle $ that produces $S$.

It only remains to see that $\left\langle B^{3},G,\epsilon \right\rangle $
has no edges of cycle $2$. This is true because, due to the construction of $%
\left\langle B^{3},G,\epsilon \right\rangle $, the $2$-components, $1$%
-components and $0$-components of $S$ coincide exactly with its faces, edges
and vertices. The interiors of the faces of $\left\langle B^{3},G,\epsilon
\right\rangle $ produce, when identified, exactly the 2-components of $S$.
Similarly, the edges of $\left\langle B^{3},G,\epsilon \right\rangle $
produce exactly the 1-components of $S$, and the vertices the 0-components.
This implies that no point in $S$\ coming from an edge of $\left\langle
B^{3},G,\epsilon \right\rangle $ has a neighborhood with the shape of a
disk, and by our analysis of neighborhoods we can conclude that $%
\left\langle B^{3},G,\epsilon \right\rangle $ has no edges of cycle $2$.

On the other hand, it is easy to see that if a polyhedron has no edges of
cycle 2, the scar produced by it is necessarily cellular. $\square $

\bigskip

\begin{lemma}
\label{L6}Let $\left\langle B^{3},G,\epsilon \right\rangle $ be a polyhedron
all whose edges are of cycle 3, all whose vertices are or order less than or
equal to 4, and all whose vertices are adjacent to exactly three edges. Then
every vertex of $\left\langle B^{3},G,\epsilon \right\rangle $ is of order
equal to 4.
\end{lemma}

\textit{Proof}. We have already proven this before unintentionally. The existence of vertices of order 1 is discarded because it
implies the existence of edges of cycle 1. Besides those vertices, if they
existed, would be needless. The existence of vertices of order 2 is
similarily descarded because it implies the existence of edges of cycle 2.
Finally, let us suppose that a vertex $v$ has order 3. Then, there
necessarily exist two identified edges of cycle 3 adjacent to $v$, and $v$
has a neighborhood with the shape of three half disks, as in Figure 2. (b)
But this implies that the three vertices in the polyhedron whose
identification turns them into $v$ are, each of them, adjacentent to exactly
two edges of cycle 3, which violates the hypotheses. We conclude then that
every vertex has order 4. $\square $

\bigskip

We will continue with the following definition. We say that a polyhedron is 
\textit{distinguished} if all of its edges are of cycle 3, all its vertices
are of order 4, and all its vertices are adjacent to exactly three edges. It
can be proved that in fact the condition of the vertices to have order 4 is
superfluous, for it is implied by the other two conditions. We are in a
position now to prove the following theorem, that is the main result of this
section, and establishes the equivalence between the scars of distinguished
polyhedra that produce 3-manifolds, and the special spines of 3-manifolds.

\begin{theorem}
\label{T3}Let $M$ be a closed, connected, orientable 3-manifold, and $%
S\subseteq M$. Then, $S$ is a special spine of $M$ if and only if it is the
scar of a distinguished polyhedron that produces $M$.
\end{theorem}

\bigskip

\textit{Proof}. It follows from Lemmas \ref{L4}, \ref{L5} and \ref{L6} that
a scar is simple and cellular if and only if there exists a distinguished
polyhedron that produces it. We will use this fact along the proof.

Let $S$ be a special spine of the manifold $M$. Let us see that there exists
a distinguished polyhedron that produces $S$ and $M$. By Theorem \ref{T2},
we know that $S$ is the scar of a polyhedron $\left\langle B^{3},H,\eta
\right\rangle $ that produces $M$. However, we have no way to know wether $%
\left\langle B^{3},H,\eta \right\rangle $ is distinguished or not, for which
this polyhedron is not of interest to us. What interest us in this regard is
the fact that the spine $S$ is also a scar. This fact, in conjunction with
the definition of special spine, implies that $S$ is furthermore a simple
and cellular scar. Therefore, there exists a distinguished polyhedron $%
\left\langle B^{3},G,\epsilon \right\rangle $ that produces $S$. Let us see
then that $\left\langle B^{3},G,\epsilon \right\rangle $ produces $M$ also.

Let $P$ be the space produced by $\left\langle B^{3},G,\epsilon
\right\rangle $. Then, since $S$ is the scar of $\left\langle
B^{3},G,\epsilon \right\rangle $ and $S$ is simple, by Lemma \ref{L3}, $P$
is a 3-manifold. Besides, by Theorem \ref{T2}, $S$ is a
spine of $P$. Moreover, since $S$ is a special spine of $M$, by its own
topology $S$ is a special spine of $P$. Now, since $S$ is a special spine of
both $P$ and $M$, we have that $P=M$; given that two manifolds with
homeomorphic special spines are necessarily homeomorphic (see \cite{Mat2}).

On the other hand, if $\left\langle B^{3},G,\epsilon \right\rangle $ is a
distinguished polyhedron that produces $M$, and if $S$ is its scar, we have
by Theorem \ref{T2} that $S$ is a spine of $M$. Moreover, since $%
\left\langle B^{3},G,\epsilon \right\rangle $ is distinguished, $S$ is a
simple and cellular scar for which it is in fact a special spine of $M$. $%
\square $

\bigskip

\section{Distinguished Polyhedra}

In this section we will prove that the distinguished polyhedra are a
presentation of the closed, connected, orientable 3-manifolds. This
presentation is in fact equivalent to the presentation by special
thickenable (PL) polyhedra, or special spines, as we will prove in Section 6
(see \cite{Mat2}).

\begin{theorem}
\label{T4}Let $\left\langle B^{3},G,\epsilon \right\rangle $ be a
distinguished polyhedron, and $S$ its scar. Then, the space produced by $%
\left\langle B^{3},G,\epsilon \right\rangle $ is a 3-manifold, and $S$ is a
special spine of such 3-manifold.
\end{theorem}

\textit{Proof}. Since $\left\langle B^{3},G,\epsilon \right\rangle $ is
distinguished, $S$ is simple and, by Lemma \ref{L3}, $\left\langle
B^{3},G,\epsilon \right\rangle $ produces a 3-manifold. By Theorem \ref{T3}, 
$S$ is a special spine of such 3-manifold. $\square $

\bigskip

\begin{theorem}
\label{T5}The distinguished polyhedra are a presentation of the closed,
connected, orientable 3-manifolds; where each polyhedron presents the
manifold that is its quotient space.
\end{theorem}

\textit{Proof}. It only remains to see that for every 3-manifold $M$ of this
type there exists a distinguished polyhedron that produces it. This is true
because $M$ has some special spine (see \cite{Mat2}), and by Theorem \ref{T3}
that spine is the scar of a distinguished polyhedron that produces $M$. $%
\square $

\bigskip

\section{Spines and Polyhedra}

Up to this point we have fully established the relation between the spines
of 3-manifolds and the scars of polyhedra. Our task now will be to establish
a more direct relation between spines and polyhedra. Specifically, we will
establish sufficient and necessary conditions for two polyhedra to produce
the same quotient space and the same scar.

Let us recall for a moment the concept of alikeness between relations given
in Definition \ref{D3}, and consider two polyhedra $\left\langle
B^{3},G,\epsilon \right\rangle $ and $\left\langle B^{3},H,\eta
\right\rangle $. For $\left\langle B^{3},G,\epsilon \right\rangle $ we
define $\check{E}q\epsilon $ as the relation obtained by adding to $%
Eq\epsilon $ every one-point set of the form $\left\{ x\right\} $, with $%
x\in \mathring{B}$. We define $\check{E}q\eta $ in the same way. Let
us notice that relations $Eq\epsilon $ and $Eq\eta $ in $\partial B^{3}$ are
alike if and only if relations $\check{E}q\epsilon $ and $\check{E}q\eta $
on $B^{3}$\ are alike.

Now, if $\left\langle B^{3},G,\epsilon \right\rangle $ and $\left\langle
B^{3},H,\eta \right\rangle $ are polyhedra for which $Eq\epsilon $ and $%
Eq\eta $ are alike, we say that $\left\langle B^{3},G,\epsilon \right\rangle 
$ and $\left\langle B^{3},H,\eta \right\rangle $\ are \textit{alike polyhedra%
}. The following lemma is clear.

\begin{lemma}
\label{L7}Alikeness between polyhedra is an equivalence relation. Besides,
alike polyhedra produce the same quotient space, and the same scar.
\end{lemma}

From now on, if $Eq\epsilon $ and $Eq\eta $ are alike we will just say that $%
Eq\epsilon =Eq\eta $. We will see now how to obtain the set of all the
polyhedra alike to a determined polyhedron $\left\langle B^{3},G,\epsilon
\right\rangle $. With that purpose we will define a move that allow us to
shift between alike polyhedra. Let us observe that given a polyhedron, and a
pair of faces $\left\{ F_{i},F_{i}^{-1}\right\} $ of such polyhedron, we can
draw a line or edge $a$\ that goes across $F_{i}$ from one side to another,
and at the same time draw a line $a^{-1}$ in $F_{i}^{-1}$ whose points are
the images of the points of $a$ under $\epsilon _{i}$, so that when $F_{i}$
and $F_{i}^{-1}$ are glued, $a$ glues with $a^{-1}$. Let us notice that this
process does not alter the polyhedron substantially, and that in this case $%
\left\{ a,a^{-1}\right\} $\ is a cycle of two edges of cycle two. We call
this process the insertion of an edge of cycle two, and we can conceive the
remotion of an edge of cycle two in a similar way. Insertion and remotion of
edges of cycle 2 will be the moves that will allow us to shift between alike
polyhedra, and formally we define them in the following way. The notation
\textquotedblleft $\,:\,$\textquotedblright\ will be used for adjacency
between cells.

Let $\left\langle B^{3},G,\epsilon \right\rangle $ be a polyhedron with
faces $\left\{ F_{1},F_{1}^{-1},...,F_{n},F_{n}^{-1}\right\} $ and
identification scheme $\epsilon =\left\{ \epsilon _{1},...,\epsilon
_{i},\epsilon _{j},...,\epsilon _{n}\right\} $. Let $\left\{
a_{1},a_{2}\right\} $ be a class of two edges of cycle 2 in $\left\langle
B^{3},G,\epsilon \right\rangle $, and $F_{i}$ and $F_{j}$ be the faces of
the polyhedron for which $F_{i}:a_{1}:F_{j}$ and $%
F_{i}^{-1}:a_{2}:F_{j}^{-1} $ holds. Let us notice that $F_{i}$ and $F_{j}$
are not necessarily different. Now, let us set $G^{\prime }=\overline{G-(%
\bar{a}_{1}\cup \bar{a}_{2})}$ and $\delta =\left\{ \epsilon
_{1},...,\epsilon _{i}\cup \epsilon _{j},...,\epsilon _{n}\right\} $. Thus,
if $G^{\prime }$ is connected, then $\left\langle B^{3},G^{\prime },\delta
\right\rangle $ is a polyhedron alike to $\left\langle B^{3},G,\epsilon
\right\rangle $, and we say that the first one is obtained from the second
one by the \textit{remotion of an edge of cycle 2}.

On the other hand, let $F_{i}$ be a face of $\left\langle B^{3},G,\epsilon
\right\rangle $. For $x\in F_{i}\cap G$ and $y\in F_{i}$, let $%
f:[0,1]\longrightarrow F_{i}$ be a continuous injective function such that $%
f(0)=x$ and $f(1)=y$. Let us define $a_{1}$ as $a_{1}=f([0,1])$, and $a_{2}$
as the image of $a_{1}$ under $\epsilon _{i}$. Then, depending on whether $y$
belongs or not to $G$, $F_{i}-a_{1}$ may have two connected components or
only one. If $F_{i}-a_{1}$ has two connected components, we denote its
closures by $F_{j}$ and $F_{k}$, with $k>j>n$. If $F_{i}-a_{1}$ has a single
connected component, we will understand that $F_{j}:=F_{k}:=F_{i}$.
Aditionally, we will define $F_{j}^{-1}$ and $F_{k}^{-1}$ as the images of $%
F_{j}$ and $F_{k}$ under $\epsilon _{i}$; and we will define $\epsilon _{j}$
and $\epsilon _{k}$ as the restrictions of $\epsilon _{i}$ to $F_{j}$ and $%
F_{k}$ respectively. Finally, let us set $G^{\prime }=G\cup a_{1}\cup a_{2}$%
, and $\delta =\left\{ \epsilon _{1},...,\epsilon _{i-1},\epsilon
_{i+1},...,\epsilon _{n},\epsilon _{j},\epsilon _{k}\right\} $. Then, under
these circumstances, $\left\langle B^{3},G^{\prime },\delta \right\rangle $
is a polyhedron alike to $\left\langle B^{3},G,\epsilon \right\rangle $, and
we say that the first one is obtained from the second one by the \textit{%
insertion of an edge of cycle 2}.

We can state now the following lemma.

\begin{lemma}
\label{L8}Two polyhedra are alike if and only if one of them can be obtained
from the other by insertion and remotion of edges of cycle 2.
\end{lemma}

\textit{Proof}. The implication from right to left is obtained directly from
the definition of insertion and remotion of edges of cycle 2. Let us see the
other implication. Let $\left\langle B^{3},G,\epsilon \right\rangle $ and $%
\left\langle B^{3},H,\eta \right\rangle $ be two polyhedra in the same
alikeness class. Then we can consider the graph $G\cup H$ in $\partial B^{3}$%
. If $G\cup H$\ is not connected, we can connect it by the insertion of an
edge of cycle 2 and obtain a new graph $J$. If $G\cup H$ is connected, we
define $J$ simply by $J=G\cup H$. Since $J$ is connected, $\left\langle
B^{3},J\right\rangle $ is a cell-divided ball. Restricting the relation $%
Eq\epsilon =Eq\eta $ to each of the faces of $\left\langle
B^{3},J\right\rangle $ we obtain relations $\delta _{1},...,\delta _{n}$,
such that $\delta $ is an identification scheme for $\left\langle
B^{3},J\right\rangle $. Hence, $\left\langle B^{3},J,\delta \right\rangle $
is a polyhedron in the alikeness class of $\left\langle B^{3},G,\epsilon
\right\rangle $ and $\left\langle B^{3},H,\eta \right\rangle $.

Let $\Gamma $ be the set of all the points $x\in \partial B^{3}$, whose
equivalence classes under $Eq\epsilon $ have cardinal different from $2$.
Let us see that $\Gamma $ is contained in $G\cap H$. From the discussion of
Section 3, about the shapes of the neighborhoods of the points in the scars,
it follows that the points in the scar of $\left\langle B^{3},G,\epsilon
\right\rangle $ coming from points in $\Gamma $ cannot have neighborhoods
(on the same scar) homeomorphic to disks, and therefore $\Gamma $ must be a
subset of $G$. Furthermore, that analysis reveals that $\Gamma $ is in fact
a subgraph of $G$. Since $Eq\epsilon =Eq\eta $, by symmetry, $\Gamma $ must
also be a subset of $H$. Thus, $\Gamma \subseteq G\cap H$.

The same argument shows that $\Gamma \subseteq G\cap J$ and that $\Gamma
\subseteq H\cap J$. Since $G$ and $H$ are subgraphs of $J$, we conclude that
both $\left\langle B^{3},G,\epsilon \right\rangle $ and $\left\langle
B^{3},H,\eta \right\rangle $ can be obtained from $\left\langle
B^{3},J,\delta \right\rangle $ by the remotion of edges of cycle 2, with
which we have proven the lemma. $\square $

\bigskip

Let $P$ be a pseudomanifold of type $\mathcal{P}_{1}o$, and let $S$ be an
homogeneous bidimensional spine of $P$. Then we say that a polyhedron
produces $(P,S)$ if the space that it produces is $P$ and its scar is $S$.
We are in a position now to prove the following theorem, that will be the
main result of this section, and establishes the equivalence between the
pairs of the form $(P,S)$ and the alikeness classes of polyhedra.

\begin{theorem}
\label{T7}Two polyhedra are alike if and only of they produce the same
quotient space and the same scar.
\end{theorem}

\textit{Proof}. We already know that alike polyhedra produce the same space
and the same scar (Lemma \ref{L7}). Let us see now that if two polyhedra
produce the same quotient space $P$,\ and the same scar $S$, then they are
alike.

The proof is based on the construction made in the proof of Theorem \ref{T2}%
, to prove that every spine is the scar of some polyhedron. Let us observe
that given a triangulation $T$ of $S$, the construction of such polyhedron,
just as it was carried out in the proof of the theorem, leads to a unique
polyhedron; for which we can denote the same by $\left\langle
B^{3},G_{T},\epsilon _{T}\right\rangle $. Let us see first that the
polyhedra obtained from $P$ and $S$ by this method starting from different
triangulations are all alike. Let us consider two triangulations $K$ and $L$
of $S$, and let us take one more triangulation, $R$, that be a common
subdivision of $K$ and $L$. Then $\left\langle B^{3},G_{R},\epsilon
_{R}\right\rangle $ can be obtained from both $\left\langle
B^{3},G_{K},\epsilon _{K}\right\rangle $ and $\left\langle
B^{3},G_{L},\epsilon _{L}\right\rangle $ by insertion of edges of cycle 2,
for which $\left\langle B^{3},G_{K},\epsilon _{K}\right\rangle $ and $%
\left\langle B^{3},G_{L},\epsilon _{L}\right\rangle $ are alike.

Now, let $\left\langle B^{3},H,\eta \right\rangle $ be an arbitrary
polyhedron that produces $M$ and $S$. Let us see that there exists a
triangulation $T$ of $S$ for which $\left\langle B^{3},H,\eta \right\rangle $
and $\left\langle B^{3},G_{T},\epsilon _{T}\right\rangle $ are alike.
Through the insertion of edges of cycle 2 in $\left\langle B^{3},H,\eta
\right\rangle $, it is possible to obtain a polyhedron $\left\langle
B^{3},H^{\prime },\eta ^{\prime }\right\rangle $ whose faces are all
triangular. Clearly $H^{\prime }$ induces a triangulation $T_{0}$ in $S$,
and it is easy to see that $\left\langle B^{3},H^{\prime },\eta ^{\prime
}\right\rangle =\left\langle B^{3},G_{T_{0}},\epsilon _{T_{0}}\right\rangle $%
. Thus, $T_{0}$ is a triangulation of $S$ such that $\left\langle
B^{3},G_{T_{0}},\epsilon _{T_{0}}\right\rangle $ and $\left\langle
B^{3},H,\eta \right\rangle $ are alike.

In this way we have that if $T$ is a triangulation of $S$, every polyhedron
that produces $M$ and $S$ is alike to $\left\langle B^{3},G_{T},\epsilon
_{T}\right\rangle $, which completes the proof. $\square $

\bigskip

The following theorem establishes the equivalence between the pairs of the
form $(M,S)$ (where $M$ is a closed, connected, orientable 3-manifold, and $%
S $ is a homogeneous bidimensional spine of $M$) and the alikeness classes
of polyhedra that produce manifolds.

\begin{theorem}
\label{C1}Let $M$ be a closed, connected, orientable 3-manifold, and suppose
that $S$ is a homogeneous bidimensional spine of $M$. Then, there exists a
unique alikeness class, such that all its polyhedra produce $(M,S)$, and
such that no other polyhedron produces $(M,S)$.
\end{theorem}

\textit{Proof}. It follows trivially from Theorems \ref{T2} and \ref{T7}. $%
\square $

\section{Special Spines and Polyhedra}

In Section 3 we established the relation between the special spines of
3-manifolds and the scars of distinguished polyhedra. In this last section,
in light of the developments of the previous section, we will aim for a more
direct relation between special spines and distinguished polyhedra. This
will be acomplished in the main result of this section, which is given in
Theorem \ref{C2}, and establishes the equivalence between the distinguished
polyhedra and the special spines of 3-manifolds (or special thickenable PL
polyhedra).

Before stating this result we shall give some definitions. Let $\left\langle
B^{3},G,\epsilon \right\rangle $ be a polyhedron, and let $\Gamma $ be
definied as in the proof of Lemma \ref{L8}, that is, as the set of all
points $x\in\partial B^{3}$, whose equivalence classes under $Eq\epsilon $
have cardinal different from $2$. Then $\Gamma $ is a subgraph of $G$, as we
saw in that proof. Since $\Gamma $ is defined from $Eq\epsilon $ and not
from $\left\langle B^{3},G,\epsilon \right\rangle $, we see that $\Gamma $
only depends on the alikeness class of $\left\langle B^{3},G,\epsilon
\right\rangle $. We call $\Gamma $ the \textit{essential graph} of the class
of $\left\langle B^{3},G,\epsilon \right\rangle $.

Let us notice that $\Gamma $ does not necessarily have to be connected or
non-empty, though it can be proven that it is empty only in the case of a
certain polyhedron for the projective space. Let us observe that if $\Gamma $
is connected and non-empty, $\left\langle B^{3},\Gamma \right\rangle $ is by
definition a cell-divided ball. Moreover, the restrictions of $Eq\epsilon $
to each of the faces of $\left\langle B^{3},\Gamma \right\rangle $ produce
an identification scheme $\gamma $ for $\left\langle B^{3},\Gamma
\right\rangle $, such that $Eq\gamma =Eq\epsilon $. Thus, if $\Gamma $ is
connected and non-empty, $\left\langle B^{3},\Gamma ,\gamma \right\rangle $
is a well defined polyhedron. We call $\left\langle B^{3},\Gamma ,\gamma
\right\rangle $ the \textit{minimum polyhedron} of the class of $%
\left\langle B^{3},G,\epsilon \right\rangle $.

Now, let $\Gamma $ be the essential graph (connected and non-empty) of some
alikeness class $A$, and let $\left\langle B^{3},G,\epsilon \right\rangle $
be a class representative for $A$. Then, the analysis of neighborhoods of
Section 3 reveals that $\Gamma $, besides being a subgraph of $G$, is the
union of the closures of all the edges of $G$ with cycle different from 2.
Hence, $\left\langle B^{3},\Gamma ,\gamma \right\rangle $ is the polyhedron
obtained by succesively removing all the edges of cycle 2 from $\left\langle
B^{3},G,\epsilon \right\rangle $. The arbitrariness of $\left\langle
B^{3},G,\epsilon \right\rangle $ implies that we have proven the following
lemma.

\begin{lemma}
\label{L9}If $\left\langle B^{3},\Gamma ,\gamma \right\rangle $ is the
minimum polyhedron of an alikeness class $A$, with $\Gamma $ connected and
non-empty, then every polyhedron in $A$ is obtained from $\left\langle
B^{3},\Gamma ,\gamma \right\rangle $ by the insertion of edges of cycle 2.
\end{lemma}

We can prove now the following theorem.

\begin{theorem}
\label{T8}Every distinguished polyhedron produces a closed, connected,
orientable 3-manifold, and its scar is a special spine of that manifold.
Inversely, if $M$ is a closed, connected, orientable 3-manifold and $S$ is a
special spine of $M$, then there exists a unique distinguished polyhedron
that produces $M$ and whose scar is $S$.
\end{theorem}

\textit{Proof}. The first statement is almost exactly the statement of
Theorem \ref{T4}. The closedness, connectedness and orientability are
implied by Theorem \ref{T1}. Let us prove now the other affirmation. Let $M$
be a (closed, connected, orientable) 3-manifold and $S$ a special spine of $%
M $. Then we have by Theorem \ref{T3} that there exists a distinguished
polyhedron $\left\langle B^{3},G,\epsilon \right\rangle $ that produces $%
(M,S)$. On the other hand we have, by Theorem \ref{C1}, that there exists a
unique alikeness class $A$ whose polyhedra produce $(M,S)$. From there it
follows that $\left\langle B^{3},G,\epsilon \right\rangle $ belongs to $A$.

The theorem will be proven if we show that $\left\langle B^{3},G,\epsilon
\right\rangle $ is the only distinguished polyhedron in $A$. Now, we know
that $\left\langle B^{3},G,\epsilon \right\rangle $, for being
distinguished, lacks edges of cycle 2. Since the essential graph of $A$ is
the union of the closures of all the edges of $G$ with cycle different from
2, we have that $G$ is the essential graph of $A$. Hence, since $%
\left\langle B^{3},G,\epsilon \right\rangle $ is a well defined polyhedron, $%
G$ is connected and non-empty, and $\left\langle B^{3},G,\epsilon
\right\rangle $ is the minimum polyhedron of $A$. From there it follows that
every polyhedron in $A$ is obtained from $\left\langle B^{3},G,\epsilon
\right\rangle $ by the insertion of edges of cycle 2 (Lemma \ref{L9}), and
that $\left\langle B^{3},G,\epsilon \right\rangle $ is the only
distinguished polyhedron in $A$. $\square $

\begin{theorem}
\label{C2}Every distinguished polyhedron has as its scar a special spine.
Inversely, for every special spine there exists a unique distinguished
polyhedron that has it as its scar.
\end{theorem}

\textit{Proof}. The first affirmation holds by Theorem \ref{T4}. Let us see
the other affirmation. Let $S$ be a special spine of some closed, connected,
orientable 3-manifold. Since $S$ is special, that 3-manifold is unique, and
we can denote it by $M_{S}$. By the previous theorem (Theorem \ref{T8}),
there exists a unique distinguished polyhedron that produces $(M_{S},S)$.
From there it follows that there exists a unique distinguished polyhedron
that produces $S$. $\square $

\bigskip

The previous theorem, which is the main result of this section, establishes
a one to one correspondence between the distinguished polyhedra and the
special spines, in which each polyhedron produces its corresponding spine.
It is a consecuence of the same theorem and Theorem \ref{T3} that, for every
distinguished polyhedron, it and its corresponding special spine present
one and the same 3-manifold. We conclude then that the presentations of
3-manifolds by distinguished polyhedra and special spines are in fact
equivalent.

\end{document}